\documentclass[11pt]{article}
\usepackage[utf8]{inputenc}
\usepackage[english]{babel}
\usepackage{amsthm}
\usepackage{graphicx}
\usepackage{url}
\usepackage{hyperref}
\usepackage{amsmath}
\usepackage{listings}

\begin{document}

\title{Analytical recurrence formulas for non-trivial zeros of the Riemann zeta function}

\author{Artur Kawalec}

\date{}
\maketitle

\begin{abstract}
In this article, we develop four types of analytical recurrence formulas for non-trivial zeros of the Riemann zeta function on critical line assuming (RH). Thus, all non-trivial zeros up to the $n$th order must be known in order to generate the $n$th+1 non-trivial zero. All the presented formulas are based on certain closed-form representations of the secondary zeta function family, which are already available in the literature. We also present a formula to generate the non-trivial zeros directly from primes. Thus all primes can be converted into an individual non-trivial zero,  and we also give a set of formulas to convert all non-trivial zeros into an individual prime. We also extend the presented results to other Dirichlet-L functions, and in particular, we develop an analytical recurrence formula for non-trivial zeros of the Dirichlet beta function. Throughout this article, we also numerically compute these formulas to high precision for various test cases and review the computed results.
\end{abstract}

\newpage
\section{Introduction}
The Riemann zeta function is classically defined by an infinite series

\begin{equation}\label{eq:20}
\zeta(s) = \sum_{n=1}^{\infty} \frac{1}{n^s}
\end{equation}
which is absolutely convergent for $\Re(s)>1$, where $s$ is a complex variable $s=\sigma+it$. By analytically continuing (1) to the whole complex plane, the function has an infinity of complex non-trivial zeros constrained to lie in a critical strip region $0<\sigma<1$. The $n$th zero is denoted as $\rho_n=\sigma_n+it_n$ and is a solution of

 \begin{equation}\label{eq:20}
\zeta(\rho_n) = 0
\end{equation}
for $n\geq 1$. The first few zeros on the critical line at $\sigma_n=\frac{1}{2}$ have imaginary components  $t_1 = 14.13472514...$, $t_2 = 21.02203964...$, $t_3 = 25.01085758...$, and so on, which were originally computed numerically using an equation solver, but if the Riemann Hypothesis (RH) is true, then they can be represented by an analytical recurrence formula as

\begin{equation}\label{eq:20}
t_{n+1} = \lim_{m\to\infty}\left[\frac{(-1)^{m}}{2}\left(2^{2m}-\frac{1}{(2m-1)!}\log (|\zeta|)^{(2m)}\big(\frac{1}{2}\big)-\frac{1}{2^{2m}}\zeta(2m,\frac{5}{4})\right)-\sum_{k=1}^{n}\frac{1}{t_{k}^{2m}}\right]^{-\frac{1}{2m}}
\end{equation}
that we developed in [7], thus all non-trivial zeros up to the $n$th order must be known in order to generate the $n$th+1 non-trivial zero. We consider the $2m$ limiting value as to ensure that it is even. The key component of this representation is a term
\begin{equation}\label{eq:20}
\log (|\zeta|)^{(2m)}\big(\frac{1}{2}\big)=\frac{d^{2m}}{ds^{2m}}\log(|\zeta(s)|)\Bigr\rvert_{s = \frac{1}{2}}
\end{equation}
which is the $2m$-th derivative of $\log[\zeta(s)]$ evaluated at $s=\frac{1}{2}$, which can also be written as a Cauchy integral as
\begin{equation}\label{eq:20}
\log (|\zeta|)^{(2m)}\big(\frac{1}{2}\big)=\frac{(2m)!}{2\pi i}\oint_{\Omega}^{}\frac{\log (|\zeta(z)|)}{(z-\frac{1}{2})^{2m+1}}dz
\end{equation}
where $\Omega$ is a small contour encircling a point $s=\frac{1}{2}$. Also, $\zeta(s,a)$ is the Hurwitz zeta function
\begin{equation}\label{eq:20}
\zeta(s,a) = \sum_{n=0}^{\infty} \frac{1}{(n+a)^s},
\end{equation}
which is a shifted version of (1) by an arbitrary parameter $a\neq 0,-1,-2,\ldots$.

The essential theme in this paper is the application of the closed-form formulas for the secondary zeta function family, which we say are a generalized zeta type of series over non-trivial zeros. What we mean by the closed-form formulas is an independent representation of such series that is not involving non-trivial zeros as was initially defined, and these formulas may be composed of a combination of known constants, different infinite series or integrals, and as we will find, the $m$th logarithmic derivatives of some function.

In no particular order, we define several variations of the secondary zeta functions. The first kind is a generalized zeta series over imaginary parts of non-trivial zeros defined by

\begin{equation}\label{eq:20}
Z_1(s) =\sum_{n=1}^{\infty} \frac{1}{t_n^s}.
\end{equation}
The second kind is a generalized zeta series over complex non-trivial zeros taken in conjugate-pairs defined as
\begin{equation}\label{eq:20}
Z_2(s) =\sum_{n=1}^{\infty} \left(\frac{1}{\rho_n^s}+\frac{1}{\bar{\rho}_n^s}\right).
\end{equation}
The value for
\begin{equation}\label{eq:20}
Z_2(1) = 1+\frac{1}{2}\gamma-\frac{1}{2}\log(4\pi)
\end{equation}
is commonly known throughout the literature [4, p.67]. The third kind is a generalized zeta series over the complex magnitude, or modulus, squared of non-trivial zeros defined as
\begin{equation}\label{eq:20}
\begin{aligned}
Z_3(s) &= \sum_{n=1}^{\infty}\frac{1}{|\rho_n|^{2s}}=\sum_{n=1}^{\infty}\frac{1}{(\frac{1}{4}+t_n^2)^{s}}.
\end{aligned}
\end{equation}
The forth kind is a generalized Jacobi theta series for the sum of exponentials over the imaginary parts of non-trivial zeros defined by

\begin{equation}\label{eq:20}
\begin{aligned}
Z_4(s) &= \sum_{n=1}^{\infty}e^{-t_n^2s}.
\end{aligned}
\end{equation}

In this article, we consider a case when $s\to\infty$ and use the closed-form representations of these family of secondary zeta functions that are already available in the literature and extract non-trivial zeros under the right excitation of these series. In the previous paper [7], we developed a formula for non-trivial zeros (3) based on the secondary zeta function of the first type (7), which we will review in Section 2. In Section 3, we develop another independent formula for non-trivial zeros given by Matsuoka [9] based on combining the secondary zeta functions for (8) and (10). In Section 4, we further expand on Section 3 and develop another representation for non-trivial zeros, which can be obtained directly from primes. Thus, we can convert all primes into an individual non-trivial zero. Furthermore, we also give a set of formulas to convert all non-trivial zeros into an individual prime. In Section 5, we outline another formula based on the Jacobi type of secondary zeta function (11). And in Section 6, we extend these formulas to the non-trivial zeros of the Dirichlet-L function, and in particular, we develop a similar formula to (3) for the non-trivial zeros of the Dirichlet beta function defined as

\begin{equation}\label{eq:1}
\beta(s) = \sum_{n=0}^{\infty}\frac{(-1)^n}{(2n+1)^s},
\end{equation}
which is convergent for $\Re(s)>0$.

The computational aspect of these closed-form formulas is itself challenging and requires very high arbitrary numerical precision. Throughout this article, we compute these formulas numerically in the PARI/GP and Mathematica software package [10][14] (and the scripts are also presented) and we discuss the computed results. We find that most of these formulas converge to non-trivial zeros to a reasonable number of digits after the decimal place. But, for the Jacobi generalized series (11), we cannot achieve convergence to non-trivial zeros as the precision required to compute it is presently outside the range of our workstation to compute directly.

\section{A formula for non-trivial zeros of the first kind}
We consider the secondary zeta function of the first type as

\begin{equation}\label{eq:20}
Z_1(s) =\sum_{n=1}^{\infty} \frac{1}{t_n^s} =\frac{1}{t_{1}^{s}}+\frac{1}{t_{2}^{s}}+\frac{1}{t_{3}^{s}}+\ldots,
\end{equation}
which is convergent for $\Re(s)>1$. The special values for the first few positive integers is

\begin{equation}\label{eq:9}
\begin{aligned}
Z_1(2) &=\frac{1}{2}(\log |\zeta|)^{(2)}\big(\frac{1}{2}\big)+\frac{1}{8}\pi^2+\beta(2)-4 \\
     &= 0.023104993115418970788933810430\dots, \\
     &\\
Z_1(3) &= 0.000729548272709704215875518569\dots, \\
     &\\
Z_1(4) &=-\frac{1}{12}(\log |\zeta|)^{(4)}\big(\frac{1}{2}\big)-\frac{1}{24}\pi^4-4\beta(4)+16 \\
     &= 0.000037172599285269686164866262\dots, \\
     &\\
Z_1(5) &= 0.000002231188699502103328640628\dots.
\end{aligned}
\end{equation}
The values for even integer argument are given by a closed-form formula for $Z_1(2m)$ as

\begin{equation}\label{eq:20}
\begin{aligned}
Z_1(2m) = (-1)^m \bigg[-\frac{1}{2(2m-1)!}(\log |\zeta|)^{(2m)}\big(\frac{1}{2}\big)+\\
         -\frac{1}{4}\left[(2^{2m}-1)\zeta(2m)+2^{2m}\beta(2m)\right]+2^{2m}\bigg]
\end{aligned}
\end{equation}
assuming (RH) and valid for a positive integer $m\geq 1$, but the $2m$ ensures that the limit variable is even. Although this formula looks relatively simple, it is a result of a complicated and detailed analysis by works of Voros [11][12] by analytically continuing (13) to the whole complex plane using Mellin transforms and tools from spectral theory. There is no known closed-form formula such as this is valid for odd integer argument. The odd values given were computed using a special algorithm that was developed by Arias De Reyna [2] (further elaborated in Appendix A) and is available as a stand-alone function in a Python library. It would otherwise take billions of zeros using (13) to reach such accuracy for the odd values.

Furthermore, we also have a useful identity

\begin{equation}\label{eq:20}
\frac{1}{2^s}\zeta\big(s,\frac{5}{4}\big)=\sum_{k=1}^{\infty}\frac{1}{\left(\frac{1}{2}+2k\right)^s}=2^s\left[\frac{1}{2}\left((1-2^{-s})\zeta(s)+\beta(s)\right)-1\right],
\end{equation}
also found in [11, p.681], in which we can express the zeta and beta terms in terms of a Hurwitz zeta function resulting in a compact form

\begin{equation}\label{eq:20}
Z_1(2m) =\frac{(-1)^{m}}{2}\left[2^{2m}-\frac{1}{(2m-1)!}\log (|\zeta|)^{(2m)}\big(\frac{1}{2}\big)-\frac{1}{2^{2m}}\zeta(2m,\frac{5}{4})\right].
\end{equation}

\noindent To find the non-trivial zeros, we consider solving for $t_1$ in (13) to obtain

\begin{equation}\label{eq:20}
\frac{1}{t_{1}^{s}}=Z_1(s) -\frac{1}{t_{2}^{s}}-\frac{1}{t_{3}^{s}}-\ldots
\end{equation}
and then we get
\begin{equation}\label{eq:20}
t_1=\left(Z_1(s) -\frac{1}{t_{2}^{s}}-\frac{1}{t_{3}^{s}}-\ldots\right)^{-1/s}.
\end{equation}
If we then consider the limit
\begin{equation}\label{eq:20}
t_1=\lim_{s\to\infty}\left(Z_1(s) -\frac{1}{t_{2}^{s}}-\frac{1}{t_{3}^{s}}-\ldots\right)^{-1/s}
\end{equation}
then, since $O[Z_1(s)]\sim O(t_1^{-s})$, and so the higher order non-trivial zeros decay as $O(t_2^{-s})$ faster than $Z_1(s)$, and so $Z_1(s)$ dominates the limit, hence we have
\begin{equation}\label{eq:20}
t_1=\lim_{s\to\infty}\left[Z_1(s)\right]^{-1/s}.
\end{equation}
Hence, by substituting (17) into (21), the formula for $t_1$ becomes

\begin{equation}\label{eq:20}
t_{1} = \lim_{m\to\infty}\left[\frac{(-1)^{m}}{2}\left(2^{2m}-\frac{1}{(2m-1)!}\log (|\zeta|)^{(2m)}\big(\frac{1}{2}\big)-\frac{1}{2^{2m}}\zeta(2m,\frac{5}{4})\right)\right]^{-\frac{1}{2m}}.
\end{equation}
Next, we numerically verify this formula in PARI, and the script is shown in Listing $1$. We broke up the representation (22) into several parts A to D. Also, sufficient memory must be allocated and precision set to high before running the script, we usually set precision to $2000$ decimal places. We utilize the Hurwitz zeta function representation since it is available in PARI and the \textbf{derivnum} function for computing the $m$th derivative very accurately for high $m$. The results are summarized in Table $1$ for various values of limit variable $m$ from low to high, where we can observe the convergence to $t_1$ as $m$ increases. Already at $m=10$ we get several digits of $t_1$, and at $m=100$ we get over $30$ digits. We performed even higher precision computations, and the result is clearly converging to $t_1$. In Appendix B, we also give a second script in the Mathematica software package to compute (22) using the Cauchy integral (5).

\begin{table}[hbt!]
\caption{The computation of $t_1$ by equation (22) for different $m$.} 
\centering 
\begin{tabular}{c c c} 
\hline\hline 
m & $t_1$ (First 30 Digits)  & Significant Digits\\ [0.5ex] 
\hline 
$1$ & 6.578805783608427637281793074245 & 0  \\
$2$ & 12.806907343833847091925940068962 & 0 \\
$3$ & 13.809741306055624728153992726341 & 0 \\
$4$ & 14.038096225961619450676758199577 & 0 \\
$5$ & 14.\underline{1}02624784431488524304946186056 & 1 \\
$6$ & 14.\underline{1}23297656314161936112154413740 & 1 \\
$7$ & 14.1\underline{3}0464459254236820197453483721 & 2 \\
$8$ & 14.1\underline{3}3083993992268169646789606564 & 2 \\
$9$ & 14.13\underline{4}077755601528384660110026302 & 3 \\
$10$ & 14.13\underline{4}465134057435907124435534843 & 3 \\
$15$ & 14.1347\underline{2}1950874675119831881762569 & 5 \\
$20$ & 14.13472\underline{5}096741738055664458081219 & 6\\
$25$ & 14.13472514\underline{1}055464326339414131271 & 9 \\
$50$ & 14.134725141734693\underline{7}89641535771021 & 16 \\
$100$ & 14.134725141734693790457251983562 & 34
\\ [1ex] 
\hline 
\end{tabular}
\label{table:nonlin} 
\end{table}

\lstset{language=C,deletekeywords={for,double},caption={PARI script for computing equation (22).},label=DescriptiveLabel,captionpos=b}
\begin{lstlisting}[frame=single]
{
    \\ set limit variable
    m = 250;

    \\ compute parameters A to D
    A =	2^(2*m);
    B = 1/factorial(2*m-1);
    C = derivnum(x=1/2,log(zeta(x)),2*m);
    D = (2^(-2*m))*zetahurwitz(2*m,5/4);

    \\ compute Z(2m)
    Z = (-1)^(m)*(1/2)*(A-B*C-D);

    \\ compute t1
    t1 = Z^(-1/(2*m));
    print(t1);
}
\end{lstlisting}

\newpage
Next we numerically compute (22) for even higher limit variable $m=250$, which yields

\begin{equation}\label{eq:20}
\begin{aligned}
t_1=14.13472514173469379045725198356247027078425711569924 & \\
     317568556746014996342980925676494901\underline{0}212214333747\ldots,
\end{aligned}
\end{equation}
which is accurate to $87$ digits.

Henceforth, in order to find the second non-trivial zero, we comeback to (13), and solving for $t_2$ yields

\begin{equation}\label{eq:20}
t_2=\lim_{s\to\infty}\left(Z_1(s) -\frac{1}{t_{1}^{s}}-\frac{1}{t_{3}^{s}}-\ldots\right)^{-1/s}
\end{equation}
and since the higher order zeros decay as $t_3^{-s}$ faster than $Z_1(s)-t_1^{-s}$, we then have

\begin{equation}\label{eq:20}
t_2=\lim_{s\to\infty}\left(Z_1(s) -\frac{1}{t_{1}^{s}}\right)^{-1/s}
\end{equation}
as $O(t_3^{-s})$ vanishes, and the zero becomes

\begin{equation}\label{eq:20}
t_{2} = \lim_{m\to\infty}\left[\frac{(-1)^{m}}{2}\left(2^{2m}-\frac{1}{(2m-1)!}\log (|\zeta|)^{(2m)}\big(\frac{1}{2}\big)-\frac{1}{2^{2m}}\zeta(2m,\frac{5}{4})\right)-\frac{1}{t_1^{2m}}\right]^{-\frac{1}{2m}},
\end{equation}
but we must know $t_1$ in advance in order to remove its contribution. A numerical computation for $m=250$ yields

\begin{equation}\label{eq:20}
\begin{aligned}
t_2=21.0220396387715549926284795938969027773\underline{3}355195796311 & \\
     4759442381621433519190301896683837161904986197676\ldots
\end{aligned}
\end{equation}
which is accurate to $38$ decimal places, and we assumed $t_1$ used was already pre-computed to $2000$ decimal places by other means. We cannot use the same $t_1$ computed earlier with the same limit variable as it will cause a self-cancellation in (26). Also, the numerical accuracy of $t_{n}$ must be much higher than $t_{n+1}$ to guarantee convergence.  And continuing on, the next zero is computed as
\begin{equation}\label{eq:20}
t_{3} = \lim_{m\to\infty}\left[\frac{(-1)^{m}}{2}\left(2^{2m}-\frac{1}{(2m-1)!}\log (|\zeta|)^{(2m)}\big(\frac{1}{2}\big)-\frac{1}{2^{2m}}\zeta(2m,\frac{5}{4})\right)-\frac{1}{t_1^{2m}}-\frac{1}{t_2^{2m}}\right]^{-\frac{1}{2m}},
\end{equation}
but we have to remove the first two zeros since $Z_1(s)-t_1^{-s}-t_2^{-s}\gg t_{4}^{-s}$. A numerical computation for $m=250$ yields

\begin{equation}\label{eq:20}
\begin{aligned}
t_3=25.010857580145688763213790992562821818659549\underline{6}5846378 & \\
     3317371101068278652101601382278277606946676481041\ldots
\end{aligned}
\end{equation}
which is accurate to $43$ decimal places, and we assumed $t_1$ and $t_2$ were used to high enough precision which was $2000$ decimal places in this example. As a result, if we define a partial secondary zeta function up to the $n$th order
\begin{equation}\label{eq:20}
Z_{1,n}(s) = \sum_{k=1}^{n}\frac{1}{t_{k}^{s}},
\end{equation}
then the $n$th+1 non-trivial zero is

\begin{equation}\label{eq:20}
t_{n+1}=\lim_{m\to\infty}\left[Z_1(m)-Z_{1,n}(m)\right]^{-1/m}
\end{equation}
because $t_{n}^{-s}\gg t_{n+1}^{-s}$, and the main formula:

\begin{equation}\label{eq:20}
t_{n+1} = \lim_{m\to\infty}\left[\frac{(-1)^{m}}{2}\left(2^{2m}-\frac{1}{(2m-1)!}\log (|\zeta|)^{(2m)}\big(\frac{1}{2}\big)-\frac{1}{2^{2m}}\zeta(2m,\frac{5}{4})\right)-\sum_{k=1}^{n}\frac{1}{t_{k}^{2m}}\right]^{-\frac{1}{2m}}.
\end{equation}

One can actually use any number of representations for $Z_1(s)$, and it will be interesting to find more efficient algorithms to compute them. And finally, we report a numerical result for $Z_1(2m)$ for $m=250$ as:
\begin{equation}\label{eq:20}
\begin{aligned}
Z_1(2m) = 7.18316934899718140841650578011166023417090863769600 & \\
     8517536818521464413577481501771580460474425539208\times 10^{-576}\ldots.
\end{aligned}
\end{equation}
From this number, we extracted the first $10$ non-trivial zeros, which are summarized in Table $2$ for $m=250$. The previous non-trivial zeros used were already known to high precision to $2000$ decimal places in order to compute the $t_{n+1}$. One cannot use the same $t_n$ obtained earlier with the same limit variable because it will cause self-cancellation, and the accuracy for $t_n$ must be much higher than $t_{n+1}$ to ensure convergence. Initially, we started with an accuracy of $87$ digits after decimal place for $t_1$, and then it dropped to $7$ to $12$ digits by the time it gets to $t_{10}$ zero. There is also a sudden drop in accuracy when the gaps get too small. Hence, these formulas are not very practical for computing higher-order zeros as large numerical precision is required, especially when we get to the first Lehmer pair at $t_{6709}=7005.06288$, the gap between next zero is about $\sim 0.04$.  Also, the average gap between zeros gets smaller as $t_{n+1}-t_{n}\sim\frac{2\pi}{\log(n)}$, making the use of this formula progressively harder and harder to compute.

\begin{table}[hbt!]
\caption{The $t_{n+1}$ computed by equation (32).} 
\centering 
\begin{tabular}{c c c c} 
\hline\hline 
$n$ & $t_{n+1}$ & $m=250$ & Significant Digits \\ [0.5ex] 
\hline 
$0$ & $t_{1}$ & 14.134725141734693790457251983562 & 87 \\
$1$ & $t_{2}$  & 21.022039638771554992628479593896 & 38 \\
$2$ & $t_{3}$ & 25.010857580145688763213790992562 & 43 \\
$3$ & $t_{4}$  & 30.424876125859513\underline{2}09940851142395 & 16 \\
$4$ & $t_{5}$  & 32.9350615877391896906623689640\underline{7}3 & 29 \\
$5$ & $t_{6}$  & 37.58617815882567125\underline{7}190902153280 & 18 \\
$6$ & $t_{7}$  & 40.918719012147\underline{4}63977678179889317 & 13 \\
$7$ & $t_{8}$  & 43.3270732809149995194961\underline{1}7449701 & 23 \\
$8$ & $t_{9}$  & 48.005150\underline{8}79831498066163921378664 & 7 \\
$9$ & $t_{10}$ & 49.77383247767\underline{2}299146155484901550 & 12
\\ [1ex] 
\hline 
\end{tabular}
\label{table:nonlin} 
\end{table}

\newpage
\section{A formula for non-trivial zeros of the second kind}
The secondary zeta function of the second kind as defined by (8) is

\begin{equation}\label{eq:20}
Z_2(s) =\sum_{n=1}^{\infty} \left(\frac{1}{\rho_n^s}+\frac{1}{\bar{\rho}_n^s}\right)=\sum_{n=1}^{\infty}\Bigg[\frac{1}{(\frac{1}{2}+it_n)^s}+\frac{1}{(\frac{1}{2}-it_n)^s}\Bigg],
\end{equation}
where the sum is taken over reciprocal complex zeros raised to power $s$ and taken in conjugate-pairs. The first few special values of this series are:

\begin{equation}\label{eq:9}
\begin{aligned}
    Z_2(1) &= 1+\frac{1}{2}\gamma-\frac{1}{2}\log(4\pi)\\
         &= 0.023095708966121033814310247906\ldots, \\
    &\\
    Z_2(2) &= 1+\gamma^2+2\gamma_1-\frac{1}{8}\pi^2\\
    &=  -0.046154317295804602757107990379\dots,\\
    &\\
    Z_2(3) &= 1+\gamma^3+3\gamma\gamma_1+\frac{3}{2}\gamma_2-\frac{7}{8}\zeta(3)\\
    &= -0.000111158231452105922762668238\dots, \\
    &\\
    Z_2(4) &= 1+\gamma^4+4\gamma^2\gamma_1+2\gamma_1^2+2\gamma\gamma_2+\frac{2}{3}\gamma_3-\frac{1}{96}\pi^4\\
    &= 0.000073627221261689518326771307\dots,\\
    &\\
    Z_2(5) &= 1+\gamma^5+5\gamma^3\gamma_1+\frac{5}{2}\gamma^2\gamma_2+\frac{5}{2}\gamma_1\gamma_2+5\gamma\gamma_1^2+\frac{5}{6}\gamma\gamma_3+\frac{5}{24}\gamma_4-\frac{31}{32}\zeta(5)\\
    &= 0.000000715093355762607735801093\dots.
\end{aligned}
\end{equation}
The value for $Z_2(1)$ is commonly known throughout the literature, and values for $Z_2(m)$ for $m>1$ also have a closed-form formula

\begin{equation}\label{eq:20}
\begin{aligned}
Z_2(m) = 1-(-1)^m2^{-m}\zeta(m)-\frac{\log(|\zeta|)^{(m)}(0)}{(m-1)!}
\end{aligned}
\end{equation}
valid for $m>1$ and is given by Voros in [12, p.73], Lehmer [8, p.23], and Matsuoka [9, p.249]. The formula is valid for even and odd index variable $m$. Another representation of (36) is given by

\begin{equation}\label{eq:20}
\begin{aligned}
Z_2(m) = 1- (1 - 2^{-m})\zeta(m) + \frac{g^c_m}{(m-1)!}
\end{aligned}
\end{equation}
where $g^c_m$ are Stieltjes Cumulants defined by Voros in [12, p.25], which are the expansion coefficients of a certain series

\begin{equation}\label{eq:20}
\log[(s-1)\zeta(s)]=-\sum_{n=1}^{\infty}\frac{(-1)^n}{n!}g_n^c(s-1)^s,
\end{equation}
as to conveniently extract the $m$th derivative at $s=0$. This series is not commonly known, but the very well-known Laurent expansion coefficients for the series
\begin{equation}\label{eq:20}
\zeta(s)=\frac{1}{s-1}+\sum_{n=0}^{\infty}\frac{(-1)^n}{n!}\gamma_n(s-1)^n
\end{equation}
are the Stieltjes constants, and the value for $\gamma_0=\gamma=0.5771256649\dots$ is the Euler-Mascheroni constant. We also define yet another series
\begin{equation}\label{eq:20}
-\frac{\zeta'(s)}{\zeta(s)}=\frac{1}{s-1}+\sum_{n=0}^{\infty}\eta_n(s-1)^n,
\end{equation}
where $\eta_n$ are its Laurent expansion coefficients. The relation between $\eta_n$ and $\gamma_n$ is given by a recurrence relation

\begin{equation}\label{eq:20}
\eta_n=(-1)^{n+1}\left[\frac{n+1}{n!}\gamma_n+\sum_{k=0}^{n-1}\frac{(-1)^{k-1}}{(n-k-1)!}\eta_k\gamma_{n-k-1}\right]
\end{equation}
found in Coffey [3, p.532] and then $g_n^c$ is
 \begin{equation}\label{eq:20}
g_n^c=(-1)^n(n-1)!\eta_{n-1}
\end{equation}
found in [12, p.25]. Essentially, all of these coefficients are variants of one another. Now, by substituting $g_n^c$ to equation (37) can now generate the values for $Z_2(m)$ in terms of the Stieltjes constants in (35) as shown on Wolfram's website [13].

Now, when trying to extract the non-trivial zeros using $Z_2(s)$, we encounter difficulty when combining the conjugate-pairs of zeros as shown

\begin{equation}\label{eq:20}
w_n(s) =\frac{1}{(\frac{1}{2}+it_n)^s}+\frac{1}{(\frac{1}{2}-it_n)^s},
\end{equation}
where $w_n(s)$ is an $n$th conjugate-pair term. From this form, it is not readily possible to separate the non-trivial zero terms in the limit as $s\to\infty$ similarly as in the non-trivial zero formula (3). The solution to this was given by Matsuoka [9], and what we can do is a slightly different interpretation. First, we find

\begin{equation}\label{eq:20}
w_n^{2}(s)=\frac{2}{(\frac{1}{4}+t_n^2)^s}+\frac{1}{(\frac{1}{2}+it_n)^{2s}}+\frac{1}{(\frac{1}{2}-it_n)^{2s}},
\end{equation}
then we get

\begin{equation}\label{eq:20}
\frac{1}{|\rho_n|^{2s}}= \frac{1}{(\frac{1}{4}+t_n^2)^{s}} = \frac{1}{2}[w_n^2(s)-w_n(2s)].
\end{equation}
This motivates to define a new secondary zeta function for the reciprocal powers of complex magnitude, or modulus, squared of $\rho$ as

\begin{equation}\label{eq:20}
\begin{aligned}
Z_3(s) &= \sum_{n=1}^{\infty}\frac{1}{|\rho_n|^{2s}}=\sum_{n=1}^{\infty}\frac{1}{(\frac{1}{4}+t_n^2)^{s}} \\
       &= \sum_{n=1}^{\infty}\frac{1}{2}[w_n^2(s)-w_n(2s)] \\
       &= \sum_{n=1}^{\infty}\frac{1}{2}w_n^2(s)-\frac{1}{2}Z_2(2s).
\end{aligned}
\end{equation}
We then need to find another formula for $w^2_n(s)$ which is an $n$th conjugate-pair term squared as defined above. If we expand (34) as
\begin{equation}\label{eq:20}
Z_2(s) = \sum_{n=1}^{\infty}w_n(s)=w_1(s)+w_2(s)+w_3(s)+\ldots \\
\end{equation}
so that
\begin{equation}\label{eq:20}
Z_2^2(s) = \sum_{n=1}^{\infty}w_n(s)=w_1^2(s)+O[2w_1(s)w_2(s)], \\
\end{equation}
because  $|w_1(s)^2|\gg|w_1(s)||w_2(s)|$ as $s\to\infty$. Now, substituting (48) to (46) yields an asymptotic formula for

\begin{equation}\label{eq:20}
Z_3(s) \sim \frac{1}{2}[Z_2^2(s)-Z_2(2s)] \\
\end{equation}
as $s\to\infty$, which is actually what we need to extract non-trivial zeros. What we don't have is a formula for $Z_3(s)$ for an arbitrary $s$, but that is not needed for the next step. Hence, if we begin with the secondary zeta function
\begin{equation}\label{eq:20}
Z_3(s) = \frac{1}{(\frac{1}{4}+t_1^2)^{s}}+\frac{1}{(\frac{1}{4}+t_2^2)^{s}}+\frac{1}{(\frac{1}{4}+t_3^2)^{s}}+\ldots
\end{equation}
and then solving for $t_1$  we obtain

\begin{equation}\label{eq:20}
 \frac{1}{(\frac{1}{4}+t_1^2)^{s}}=Z_3(s)-\frac{1}{(\frac{1}{4}+t_2^2)^{s}}-\frac{1}{(\frac{1}{4}+t_3^2)^{s}}-\ldots
\end{equation}
and then we get
\begin{equation}\label{eq:20}
\frac{1}{4}+t_1^2=\left[Z_3(s) -\frac{1}{(\frac{1}{4}+t_2^2)^{s}}-\frac{1}{(\frac{1}{4}+t_3^2)^{s}}-\ldots\right]^{-1/s}
\end{equation}
and this leads to
\begin{equation}\label{eq:20}
t_1=\left[\left(Z_3(s) -\frac{1}{(\frac{1}{4}+t_2^2)^{s}}-\frac{1}{(\frac{1}{4}+t_3^2)^{s}}-\ldots\right)^{-1/s}-\frac{1}{4}\right]^{1/2}.
\end{equation}
If we then consider the limit as $s\to\infty$, then the higher order terms decay as $O[(\frac{1}{4}+t^2_2)^{-s}]$, and hence, substituting equation (49) for $Z_3(s)$ yields
\begin{equation}\label{eq:20}
t_1=\lim_{s\to\infty}\left[\left(\frac{1}{2}Z_2^2(s)-\frac{1}{2}Z_2(2s)\right)^{-1/s}-\frac{1}{4}\right]^{1/2}
\end{equation}
which was given by Matsuoka in [9]. One can substitute any representation of $Z_2(s)$ such as by equations (36) or (37), which involves the Stieltjes constants expansion as shown in Appendix C for the first few $m$.

Next, we numerically verify this formula in PARI, and the script is shown in Listing $2$. We use equation (36) for $Z_2(s)$ and broke up the representation (54) into several parts A to C. And as before, sufficient memory must be allocated and precision set to high before running the script. The results are summarized in Table $3$ for various limit values of $m$ from low to high, and we can observe the convergence to the real value as $m$ increases. Already at $m=10$ we get several digits of $t_1$, and at $m=100$ we get over $16$ digits. We observe that for odd value of $m$ the convergence is slightly better than for even $m$. We performed higher precision computations, and the result is clearly converging to $t_1$.

Next we numerically verify it for $m=250$ which yields
\begin{equation}\label{eq:20}
t_1= 14.134725141734693790457251983562470270784257\underline{1}2346050\ldots
\end{equation}
which is accurate to $43$ digits. Henceforth, in order to find the second non-trivial zero, we comeback to (50), and solving for $t_2$ yields

\begin{equation}\label{eq:20}
t_2=\left[\left(Z_3(s) -\frac{1}{(\frac{1}{4}+t_1^2)^{s}}-\frac{1}{(\frac{1}{4}+t_3^2)^{s}}\ldots\right)^{-1/s}-\frac{1}{4}\right]^{1/2}
\end{equation}
and since the higher order zero terms decay as $\frac{1}{(\frac{1}{4}+t_3^2)^{s}}$ faster than $Z_3(s)-\frac{1}{(\frac{1}{4}+t_1^2)^{s}}$, we then have

\begin{equation}\label{eq:20}
t_2=\lim_{s\to\infty}\left[\left(\frac{1}{2}Z_2^2(s)-\frac{1}{2}Z_2(s)-\frac{1}{(\frac{1}{4}+t_1^2)^{s}}\right)^{-1/s}-\frac{1}{4}\right]^{1/2}.
\end{equation}
And continuing on, the next zero is computed as

\begin{equation}\label{eq:20}
t_3=\lim_{s\to\infty}\left[\left(\frac{1}{2}Z_2^2(s)-\frac{1}{2}Z_2(2s)-\frac{1}{(\frac{1}{4}+t_1^2)^{s}}-\frac{1}{(\frac{1}{4}+t_2^2)^{s}}\right)^{-1/s}-\frac{1}{4}\right]^{1/2}
\end{equation}
since the higher order zero terms decay as $\frac{1}{(\frac{1}{4}+t_4^2)^{s}}$. As a result, if we define a partial secondary zeta function up to the $n$th order
\begin{equation}\label{eq:20}
Z_{3,n}(s) = \sum_{k=1}^{n}\frac{1}{|\rho_k|^{2s}}=\sum_{k=1}^{n}\frac{1}{(\frac{1}{4}+t_k^2)^{s}}
\end{equation}

\noindent then the $n$th+1 non-trivial zero is

\begin{equation}\label{eq:20}
t_{n+1}=\lim_{s\to\infty}\left[Z_3(s)-Z_{3,n}(s)\right]^{-1/s}
\end{equation}
and the main recurrence formula:

\begin{equation}\label{eq:20}
t_{n+1}=\lim_{s\to\infty}\left[\left(\frac{1}{2}Z_2^2(s)-\frac{1}{2}Z_2(2s)-\sum_{k=1}^{n}\frac{1}{(\frac{1}{4}+t_k^2)^{s}}\right)^{-1/s}-\frac{1}{4}\right]^{1/2}.
\end{equation}

\newpage

\begin{table}[hbt!]
\caption{The computation of $t_1$ by equation (54) for different $m$.} 
\centering 
\begin{tabular}{c c c} 
\hline\hline 
m & $t_1$ (First 30 Digits)  & Significant Digits\\ [0.5ex] 
\hline 
$2$ & 5.561891787634141032446012810136 & 0 \\
$3$ & 13.757670503723662711511861003244 & 0 \\
$4$ & 12.161258748655529488677538477512 & 0 \\
$5$ & 14.075935317783371421926582853327 & 0 \\
$6$ & 13.579175424560852302300158195372 & 0 \\
$7$ & 14.\underline{1}16625853057249358432588137893 & 1 \\
$8$ & 13.961182494234115467191058505224 & 0 \\
$9$ & 14.\underline{1}26913415083941105873032355837 & 1 \\
$10$ &14.077114859427980275510456957007 & 0 \\
$15$ &14.1\underline{3}3795710050725394699252528681 & 2 \\
$20$ &14.13\underline{4}370485636531946259958638820 & 3 \\
$25$ &14.134\underline{7}00629574414322701677282886 & 4 \\
$50$ &14.13472514\underline{1}835685792188021492482 & 9 \\
$100$&14.134725141734693\underline{7}89329888107217 & 16
\\ [1ex] 
\hline 
\end{tabular}
\label{table:nonlin} 
\end{table}

Next, when we attempt to numerically verify (61) for higher zeros starting with a limit variable $m=250$, then we get $t_1$ accurate to 43 decimal places as before. However, such precision is not enough to compute $t_2$, so we have to increase the limit variable $m$ to achieve higher accuracy, which presently is at the limit of our test computer. We did, however, verify (61) successfully by pre-computing $Z_3(s)$ using 100 non-trivial zeros known to high precision (2000 decimal places). Then we computed the next zeros by (61), but presently, limitations of the test computer prevent computing $Z_3(s)$ using (49) to high enough precision.

\newpage

\lstset{language=C,deletekeywords={for,double},caption={PARI script for computing equation (54).},label=DescriptiveLabel,captionpos=b}
\begin{lstlisting}[frame=single]
{
    \\ set limit variable
    m1 = 250;

    \\ compute parameters A1 to C1 for Z1
    A1 =  derivnum(x=0,log(zeta(x)),m1);
    B1 = 1/factorial(m1-1);
    C1 = 1-(-1)^m1*2^(-m1)*zeta(m1);
    Z1 = C1-A1*B1;

    \\ compute 2m limit variable
    m2 = 2*m1;

    \\ compute parameters A2 to C2 for Z2
    A2 =  derivnum(x=0,log(zeta(x)),m2);
    B2 = 1/factorial(m2-1);
    C2 = 1-(-1)^m2*2^(-m2)*zeta(m2);
    Z2 = C2-A2*B2;

    \\ compute t1
    t1 = (((Z1^2-Z2)/2)^(-1/m1)-1/4)^(1/2);
    print(t1);
}
\end{lstlisting}

\section{Non-trivial zeros from primes}
In this section, we develop a variation of a formula for non-trivial zeros based on primes. We define a (Hurwitz) shifted version of $Z_2(s)$ by a parameter $a$ as

\begin{equation}\label{eq:20}
\begin{aligned}
Z_2(s|a) &=\sum_{n=1}^{\infty} \left[\frac{1}{(\rho_n+a-\frac{1}{2})^s}+\frac{1}{(\bar{\rho}_n+a-\frac{1}{2})^s}\right]\\
         &=\sum_{n=1}^{\infty} \left[\frac{1}{(a+it_n)^s}+\frac{1}{(a-it_n)^s}\right].
\end{aligned}
\end{equation}
The usual $Z_2(s)$ as defined by equation (8) is a special case for $a=\frac{1}{2}$. But when $a>\frac{1}{2}$, there is another closed-form representation

\begin{equation}\label{eq:1}
Z_2(s|a)=(a-\frac{1}{2})^{-s}-2^{-s}\zeta(s,\frac{5}{4}+\frac{1}{2}a)-\frac{1}{\Gamma(s)}\sum_{k=2}^{\infty}\frac{\Lambda(k)}{k^{\frac{1}{2}+a}}(\log k)^{s-1}
\end{equation}
found in Voros [12, p.56] which involves the von Mangoldt's function:

\begin{equation}
\Lambda(n)= \left \{
\begin{aligned}
&\log p, &&\text{if}\ n=p^k \text{ for some prime and integer } k\geq 1, \\
&0 && \text{otherwise}.
\end{aligned} \right.
\end{equation}
Now, if we apply the same arguments as in Section 3 to extract the non-trivial zeros, we obtain the first zero

\begin{equation}\label{eq:20}
t_{1}=\lim_{s\to\infty}\left[\left(\frac{1}{2}Z_2^2(s|a)-\frac{1}{2}Z_2(2s|a)\right)^{-1/s}-a^2\right]^{1/2}
\end{equation}
and the full recurrence formula:

\begin{equation}\label{eq:20}
t_{n+1}=\lim_{s\to\infty}\left[\left(\frac{1}{2}Z_2^2(s|a)-\frac{1}{2}Z_2(2s|a)-\sum_{k=1}^{n}\frac{1}{(a^2+t_k^2)^{s}}\right)^{-1/s}-a^2\right]^{1/2}.
\end{equation}
Hence through these formulas, the primes are directly converted into non-trivial zeros by an infinite series involving the $\Lambda(k)$ and the Hurwitz zeta function. This formula is valid for an arbitrary parameter $a>\frac{1}{2}$, but we find numerically that the convergence is very slow due to the nature of the von Mangoldt's function series, which requires billions of primes to reach some reasonable accuracy. When we test this formula, we find that convergence is improved when $a$ is increased but not too much in relation to the limit variable $s$. The script in PARI is shown in Listing 3, and we run equation (65) with parameters $k=10^9$ (up to a billionth value for $\Lambda(k)$) and $s=50$ and $a=15$.  The result is:

\begin{equation}\label{eq:20}
t_{1}=14.134\underline{7}3892414862370135\ldots
\end{equation}
which is accurate to 4 digits.

As shown in [7], we also outline the duality between primes and non-trivial zeros. The formula (66) converts all primes into an individual non-trivial zero. To complete the duality, it is also possible to convert all non-trivial zeros into an individual prime using the Golomb's formula for primes and the Hadamard product for $\zeta(s)$. Let $p_1=2$, $p_2=3$, $p_3=5$ and so on, define a sequence of primes, and $Q_{n}(s)$ define a partial Euler prime product up the $n$th order

\begin{equation}\label{eq:20}
Q_{n}(s)=\prod_{k=1}^{n}\left(1-\frac{1}{p_k^s}\right)^{-1}
\end{equation}
for $n>1$ and $Q_0(s)=1$, then the Golomb's recurrence formula for the $p_{n+1}$ prime is
\begin{equation}\label{eq:20}
p_{n+1}=\lim_{s\to \infty}\left[\zeta(s)-Q_n(s)\right]^{-1/s}
\end{equation}
as shown in [5] and [6]. And since $\zeta(s)$ can be written in terms of the Hadamard product in terms of non-trivial zeros
\begin{equation}\label{eq:20}
\zeta(s)=\frac{\pi^{s/2}}{2(s-1)\Gamma(1+\frac{s}{2})}\prod_{\rho}^{}\left(1-\frac{s}{\rho}\right),
\end{equation}
we can substitute (70) to (69) and obtain

\begin{equation}\label{eq:20}
p_{n+1}=\lim_{s\to \infty}\left[{\frac{\pi^{s/2}}{2(s-1)\Gamma(1+\frac{s}{2})}\prod_{\rho}^{}\left(1-\frac{s}{\rho}\right)}-Q_n(s)\right]^{-1/s}
\end{equation}
and the full recurrence formula:

\begin{equation}\label{eq:20}
p_{n+1}=\lim_{s\to \infty}\left[{\frac{\pi^{s/2}}{2(s-1)\Gamma(1+\frac{s}{2})}\prod_{\rho}^{}\left(1-\frac{s}{\rho}\right)}-\prod_{k=1}^{n}\left(1-\frac{1}{p_k^s}\right)^{-1}\right]^{-1/s}.
\end{equation}
Hence, this is a way for converting non-trivial zeros to the primes and without assuming (RH), as the Hadamard product is taken over all zeros and in conjugate-pairs.

\newpage
\lstset{language=C,deletekeywords={for,if,double},caption={PARI script for computing equation (65).},label=DescriptiveLabel,captionpos=b}
\begin{lstlisting}[frame=single]
\\ Define von Mangoldt function
Mangoldt(n)=
{
  ispower(n,,&n);
  if(isprime(n),log(n),0)
}

\\ main
{
  \\ set variables
  s1 = 50;        \\ limit variable
  a=15;           \\ arbitrary parameter
  k = 1000000000; \\ von Mangoldt sum limit

  \\ compute parameters A to C for Z1
  y1 = sum(n=2,k,Mangoldt(n)/n^(1/2+a)*log(n)^(s1-1));
  Z1 = (a-1/2)^(-s1)-2^(-s1)*
        zetahurwitz(s1,5/4+1/2*a)-1/gamma(s1)*y1;

  \\ compute double limit variable
  s2 = 2*s1;

  \\ compute parameters A to C for Z2
  y2 = sum(n=2,k,Mangoldt(n)/n^(1/2+a)*log(n)^(s2-1));
  Z2 = (a-1/2)^(-s2)-2^(-s2)*zetahurwitz(s2,5/4+1/2*a)
    -1/gamma(s2)*y2;

  \\ compute t1
  t1 = (((Z1^2-Z2)/2)^(-1/s1)-a^2)^(1/2);
  print(t1);
}
\end{lstlisting}

\section{Non-trivial zeros from infinite series over exponentials}
In this section, we explore yet another formula for non-trivial zeros. The Jacobi generalized series over the exponentials of $t_n$ is defined by

\begin{equation}\label{eq:20}
Z_4(s) = \sum_{n=1}^{\infty}e^{-t_n^2s}=e^{-t_1^2s}+e^{-t_2^2s}+e^{-t_3^2s}\ldots,
\end{equation}
which has a closed-form representation given by

\begin{equation}\label{eq:20}
Z_4(s) = A(s)-B(s),
\end{equation}
where
\begin{equation}\label{eq:20}
A(s)=-\frac{1}{2\sqrt{\pi s}}\sum_{k=2}^{\infty}\frac{\Lambda(k)}{\sqrt{k}}e^{-\frac{1}{4s}\log^2 k}+e^{\frac{s}{4}}
\end{equation}
and
\begin{equation}\label{eq:20}
B(s)=\frac{\gamma+\log (16\pi^2s)}{8\sqrt{\pi s}}-\frac{1}{4\sqrt{\pi s}}\int_{0}^{\infty}e^{-\frac{u^2}{16s}}\left(\frac{1}{u}-\frac{e^{\frac{3}{4}{u}}}{e^{u}-1}\right)du
\end{equation}
which is given in [2, p.3]. It is seen that it also involves the von Mangoldt's function and hence the primes. The terms of this series decay extremely fast due to the exponential nature. The first term is
\begin{equation}\label{eq:20}
Z_4(s) \sim O(e^{-t_1^2s})
\end{equation}
so that is suffices to solve for $t_1$ and we get

\begin{equation}\label{eq:20}
t_1 = \lim_{s\to\infty}\sqrt{-\frac{1}{s}\log Z_4(s)}
\end{equation}
and the recurrence formula is

\begin{equation}\label{eq:20}
t_{n+1} = \lim_{s\to\infty}\sqrt{-\frac{1}{s}\log\left(Z_4(s)-\sum_{k=1}^{n}e^{-t_k^2s}\right)}.
\end{equation}

Next, we numerically compute equation (73) for $s=2$ by summing the first two zeros to obtain

\begin{equation}\label{eq:20}
Z_4(s) =  2.912164200241304158784992817748\times 10^{-174}\ldots.
\end{equation}
The result converges to an extremely small value as $s$ increases, and hence the first term involving $t_1$ dominates the series. Now, if we re-compute it again using equation (74) (the script is not given) with $\Lambda(n)$ summed to $k=10^7$, then we get for

\begin{equation}\label{eq:20}
A = 0.3946415860608135898036962860711\ldots
\end{equation}
and
\begin{equation}\label{eq:20}
B = 0.394641583198706998425270589196\ldots,
\end{equation}
and then difference results in
\begin{equation}\label{eq:20}
Z_4(s) = A-B =  2.862106591378425696874573151789\times 10^{-9}\ldots.
\end{equation}
We observe that the difference here between $A$ and $B$ is on the order of $10^{-9}$ which is far too small to extract $t_1$ which is on the order of $10^{-174}$. The difference between $A$ and $B$ must be occurring very far out in the decimal places in order to converge to (80), which is presently outside of the reach of present numerical algorithms used. Also, the $\Lambda(n)$ series in $A$ is very slow to converge, while the integral term in $B$ is much faster. Hence this formula is not practical and is presently outside the range of what we can compute, but in principle, it should yield the non-trivial zeros.

\section{Non-trivial zeros of Dirichlet beta function}
The Dirichlet beta functions as defined as

\begin{equation}\label{eq:1}
\beta(s) = \sum_{n=0}^{\infty}\frac{(-1)^n}{(2n+1)^s},
\end{equation}
which is convergent for $\Re(s)>0$. It is useful to define

\begin{equation}\label{eq:1}
\beta(s) = \frac{1}{4^s}[\zeta(s,\frac{1}{4})-\zeta(s,\frac{3}{4})]
\end{equation}
in terms of the Hurwitz zeta function since it is available in most mathematical software packages where it can be efficiently computed, except at $s=1$ where it has a pole, but it could be handled in a limiting sense $s\to 1$. The value for $\beta(1)=\frac{\pi}{4}$.

Let $\rho_B=\frac{1}{2}+ir_n$ be non-trivial zeros of $\beta(s)$ on the critical line. The first few non-trivial zeros on the critical line have imaginary components  $r_1 = 6.02094890...$, $r_2 = 10.24377030...$, $r_3 = 12.98809801...$ which were originally computed numerically, but now can also be computed analytically by essentially the same arguments as described in Section 2. If we define the secondary beta function
\begin{equation}\label{eq:20}
B(s) =\sum_{k=1}^{\infty} \frac{1}{r_k^s}
\end{equation}
so that $B(s)$ is a sum of reciprocal powers of imaginary components of non-trivial zeros. Then, we consider a partial secondary beta function up to the $n$th order
\begin{equation}\label{eq:20}
B_n(s) =\sum_{k=1}^{n} \frac{1}{r_k^s},
\end{equation}
and because $r_1^{-s}\gg r_2^{-s}\gg r_3^{-s}\gg r_n^{-s}\ldots$ as $s\to\infty$, then the non-trivial zeros are given by a recursive relationship
\begin{equation}\label{eq:20}
r_{n+1} = \lim_{s\to\infty}[B(s)-B_n(s)]^{-1/s}.
\end{equation}
It now suffices to find a suitable formula for $B(s)$ which is also given by Voros in [12, p.110, Tab 10.3] as a general formula for Dirichlet-L functions. If we take $L_{\chi}$ to be $\beta(s)$ then we have

\begin{equation}\label{eq:20}
B(2m) = -\frac{1}{2}[(2^{2m}-1)\zeta(2m)+(1-2a)2^{2m}\beta(2m)]-\frac{1}{(2m-1)!}\log(|\beta|)^{(2m)}(\frac{1}{2})
\end{equation}
assuming (GRH) for $\beta(s)$. The parity parameter $a$ is related to the Dirichlet character, which we take it to be $a=1$. We further obtain

\begin{equation}\label{eq:20}
B(2m) = \frac{(-1)^{m+1}}{2}\left[2^{2m-1}[(1-2^{-2m})\zeta(2m)-\beta(2m)]+\frac{1}{(2m-1)!}\log(|\beta|)^{(2m)}(\frac{1}{2})\right],
\end{equation}
but there was probably a missing factor $(-1)^m/2$ in the original formula (89), so we reinserted it here. This form (90) is already good as is, but we proceed with some additional simplifications. Since the zeta term in (90) above is related to Dirichlet lambda function
\begin{equation}\label{eq:20}
(1-2^{-2m})\zeta(2m) = \lambda(2m) = \sum_{n=0}^{\infty}\frac{1}{(2n+1)^{2m}}
\end{equation}
we can simplify this further, and obtain
\begin{equation}\label{eq:20}
\lambda(2m)-\beta(2m) = 2\sum_{n=1}^{\infty}\frac{1}{(4n-1)^{2m}}=\frac{1}{2^{4m-1}}\zeta(2m,\frac{3}{4})
\end{equation}
which leads to a more compact form
\begin{equation}\label{eq:20}
B(2m) = \frac{(-1)^{m+1}}{2}\left[\frac{1}{(2m-1)!}\log(|\beta|)^{(2m)}(\frac{1}{2})+\frac{1}{2^{2m}}\zeta(2m,\frac{3}{4})\right].
\end{equation}
This results in a direct formula for $r_1$ as

\begin{equation}\label{eq:20}
r_1 = \lim_{m\to\infty}\left[\frac{(-1)^{m+1}}{2}\left(\frac{1}{(2m-1)!}\log(|\beta|)^{(2m)}(\frac{1}{2})+\frac{1}{2^{2m}}\zeta(2m,\frac{3}{4})\right)\right]^{-\frac{1}{2m}}
\end{equation}
and a full recurrence formula:

\begin{equation}\label{eq:20}
r_{n+1} = \lim_{m\to\infty}\left[\frac{(-1)^{m+1}}{2}\left(\frac{1}{(2m-1)!}\log(|\beta|)^{(2m)}(\frac{1}{2})+\frac{1}{2^{2m}}\zeta(2m,\frac{3}{4})\right)-\sum_{k=1}^{n} \frac{1}{r_k^{2m}}\right]^{-\frac{1}{2m}}.
\end{equation}

\begin{table}[hbt!]
\caption{The computation of $r_1$ by equation (94) for different $m$.} 
\centering 
\begin{tabular}{c c c} 
\hline\hline 
m & $r_1$ (First 30 Digits)  & Significant Digits\\ [0.5ex] 
\hline 
$1$ &  3.580234150633150009323781248620 & 0\\
$2$ &  5.728146231328241287248341234017 & 0 \\
$3$ &  5.966325900475084327722380500980 & 0 \\
$4$ &  6.\underline{0}08324723727322086185645916842 & 1 \\
$5$ &  6.\underline{0}17679912591888584424309703505 & 1 \\
$6$ &  6.02\underline{0}043240987781794077733596855 & 3 \\
$7$ &  6.02\underline{0}686849217175746999931646806 & 3 \\
$8$ &  6.02\underline{0}870797143727883542664755767 & 3 \\
$9$ &  6.020\underline{9}25125780393360202282926513 & 4 \\
$10$ & 6.0209\underline{4}1550676489284027261163265 & 5 \\
$15$ & 6.02094\underline{8}880761787735639621551287 & 6 \\
$20$ & 6.020948904\underline{6}09320778839216887766 & 10 \\
$25$ & 6.02094890469\underline{7}249155966074566560 & 12 \\
$50$ & 6.02094890469759665490251\underline{1}020221 & 24 \\
$100$ & 6.020948904697596654902511521612 & 47
\\ [1ex] 
\hline 
\end{tabular}
\label{table:nonlin} 
\end{table}

\newpage
A script to compute $r_1$ is presented in Listing 4, and the calculated values for various limit values of $m$ from low to high are shown in Table 4, where we observe a convergence to $r_1$. We also performed a very high precision computation of this formula, and the result is clearly converging to the zero.  We can also recursively compute the next higher order zeros, but as before, such a numerical computation is becoming even more difficult. But nevertheless, these formulas are indeed a closed-form representations for the zeros.

\lstset{language=C,deletekeywords={for,double},caption={PARI script for computing equation (94).},label=DescriptiveLabel,captionpos=b}
\begin{lstlisting}[frame=single]
\\ Define Dirichlet beta function
beta(x)=
{
    4^-x*(zetahurwitz(x,1/4)-zetahurwitz(x,3/4));
}

{
   \\ set limit variable
    m = 250;

    \\ compute parameters A to D
    A = derivnum(x=1/2,log(beta(x)),2*m);	
    B = 1/factorial(2*m-1);
    C = 2^(-2*m)*zetahurwitz(2*m,3/4);

    \\ compute B(2m)
    Z = (-1)^(m+1)/2*(A*B+C);

    \\ compute r1
    r1 = Z^(-1/(2*m));
    print(r1);
}

\end{lstlisting}

\texttt{Email: art.kawalec@gmail.com}

\newpage
\section{Appendix A}
There is an available algorithm developed by Arias De Reyna as described in [2] to compute $Z_1(s)$ and which is fully implemented in a Python library \textbf{mpmath} as a secondary zeta function and which is also analytically extended to the whole complex plane. Roughly, the algorithm takes a finite number of non-trivial zeros and a finite number of prime terms for the von Mangoldt's function term and estimates the remainder yielding an accurate computation of $Z_1(s)$ to a high number of decimal places. We tested this function and computed the odd values for $Z_1(s)$ accurately as shown in (14), which otherwise would take billions of zeros to reach. Listing 5 shows the Python script, and in Table 5, we compare the results for even values for $Z_1(s)$ with the reference values computed by the closed-form formula (15) as $Z_{1,ref}$ in PARI with precision set to $2000$ decimal places and give the values to the first $30$ decimal places. Also, the precision in Python was set to $100$ decimal places. The remainder output of the Python script was compared to $|Z_{1,ref}-Z_{1,python}|$.  We observe that the remainder for even values as given by the Python script matches the values computed by the closed-form formula (15), and in fact, exceeds in all the cases by at least an order of magnitude. Therefore we conclude that the given odd values should also be accurate within that remainder. In fact, $Z_1(s)$ can be computed for any complex $s$, which is further explored in [2].
\lstset{language=Python,caption={Python script for computing $Z_1(s)$ by the algorithm in \text{[2]}.},label=DescriptiveLabel,captionpos=b}
\begin{lstlisting}[frame=single]
from mpmath import *
mp.pretty = True; mp.dps = 100

z = secondzeta(3, error=True)
\end{lstlisting}

\begin{table}[hbt!]
\caption{The first 30 digits $Z_{1}(m)$ computed using the Python script with default parameter $a=0.0015$.} 
\centering 
\begin{tabular}{c c c c} 
\hline\hline 
$m$ & $Z_{1,python}(m)$  & Remainder & $|Z_{1,ref}-Z_{1,python}|$\\ [0.5ex] 
\hline 
$2$ & 0.023104993115418970788933810430 & $ 10^{-76}$ & $ 10^{-77}$ \\
$3$ & 0.000729548272709704215875518569 & $ 10^{-76}$ & - \\
$4$ & 0.000037172599285269686164866262 & $ 10^{-77}$ & $ 10^{-79}$ \\
$5$ & 0.000002231188699502103328640628 & $ 10^{-78}$ & - \\
$6$ & 0.000000144173931400973279695381 & $ 10^{-80}$ & $ 10^{-81}$ \\
$7$ & 0.000000009675344542702350408719 & $ 10^{-81}$ & - \\
$8$ & 0.000000000663031680252990869873 & $ 10^{-82}$ & $ 10^{-83}$ \\
$9$ & 0.000000000045991912392894862969 & $ 10^{-83}$ & - \\
$10$ & 0.00000000000321366415061660121 & $ 10^{-84}$ & $ 10^{-86}$ \\
$11$ & 0.00000000000022556506251559664 & $ 10^{-86}$ & -
\\ [1ex] 
\hline 
\end{tabular}
\label{table:nonlin} 
\end{table}

 \newpage
\section{Appendix B}
In this section, we compute $t_1$ in Mathematica using $Z_1(s)$ by equation (3), but instead of computing the $m$th derivative, we compute the Cauchy integral (5) taken along a closed contour $\Omega$ which is a square loop with $0.5$ sides around a point $s=\frac{1}{2}$. The script is shown in Listing 6. With $m=50$ we obtain

\begin{equation}\label{eq:20}
t_1=14.134725141734693\underline{7}896415358\ldots
\end{equation}
which is accurate to $16$ decimal places.

\newpage
\lstset{language=Mathematica,deletekeywords={WorkingPrecision,Factorial,K,Pi,I},caption={Mathematica script for computing $t_1$ by equations (3) and (5)},label=DescriptiveLabel,captionpos=b}
\begin{lstlisting}[frame=single]
(* Set limit variable *)
m = 50;

(* Define integrand *)
F[z_] := Log[Zeta[z]]/(z - 1/2)^(2 m + 1);

(* Integrate around a closed loop *)
I1 = NIntegrate[F[z],{z, 0.75 - 0.25 I, 0.75 + 0.25 I},
   WorkingPrecision -> 200];
I2 = NIntegrate[F[z],{z, 0.75 + 0.25 I, 0.25 + 0.25 I},
   WorkingPrecision -> 200];
I3 = NIntegrate[F[z],{z, 0.25 + 0.25 I, 0.25 - 0.25 I},
   WorkingPrecision -> 200];
I4 = NIntegrate[F[z],{z, 0.25 - 0.25 I, 0.75 - 0.25 I},
   WorkingPrecision -> 200];

K = Factorial[2 m]/(2 \[Pi] I);
Ix = K (I1 + I2 + I3 + I4);

(* Compute Z *)
A = Ix/Factorial[2 m - 1];
B = (2^(-2 m)) HurwitzZeta[2 m, 5/4];
Z = (-1)^(m) (1/2) (2^(2 m) - A - B);

(* Compute t1 *)
t1 = Z^(-1/(2 m))
\end{lstlisting}

\newpage
\section{Appendix C}
In this section, we expand the formulas for $t_1$ in terms of Stieltjes constants using equations (37), (41), (42) and (54). We utilize the Mathematica software package to expand the terms.
For $m=2$ be obtain an expansion:

\begin{equation}\label{eq:20}
\begin{aligned}
t_1 &\approx \left[\left(2\gamma_1-\frac{\pi^2\gamma_1}{4}+\gamma_1^2-\gamma\gamma_2-\frac{\gamma_3}{3}+\gamma^2-\frac{\pi^2}{8}-\frac{\gamma^2\pi^2}{8}+\frac{5\pi^4}{384}\right)^{-\frac{1}{2}}-\frac{1}{4}\right]^{\frac{1}{2}}
\end{aligned}
\end{equation}
\begin{equation}\label{eq:20}
t_1\approx 5.561891787634141032446012810136\ldots.\nonumber
\end{equation}
For $m=3$ we obtain an expansion:

\begin{equation}\label{eq:20}
\begin{aligned}
t_1 \approx & \Bigg[\Bigg(\gamma^3-\frac{21}{8}\gamma\gamma_1\zeta(3)-\frac{21}{16}\gamma_2\zeta(3)+3\gamma\gamma_1-\gamma_1^3+\frac{3}{2}\gamma_2+\frac{3}{2}\gamma\gamma_1\gamma_2+\\
&+\frac{3}{4}\gamma_2^2-\frac{1}{2}\gamma^2\gamma_3-\frac{1}{2}\gamma_1\gamma_3-\frac{1}{8}\gamma\gamma_4-\frac{1}{40}\gamma_5-\frac{7}{8}\zeta(3)-\frac{7}{8}\gamma^3\zeta(3)+\\
&+\frac{49}{128}\zeta(3)^2+\frac{1}{1920}\pi^6\Bigg)^{-\frac{1}{3}}-\frac{1}{4}\Bigg]^{\frac{1}{2}}
\end{aligned}
\end{equation}

\begin{equation}\label{eq:20}
t_1\approx 13.757670503723662711511861003244\ldots.\nonumber
\end{equation}
For $m=4$ we obtain an expansion:

\begin{equation}\label{eq:20}
\begin{aligned}
t_1 \approx & \Bigg[\Bigg(4\gamma^2\gamma_1-\frac{1}{24}\gamma^2\gamma_1\pi^4+2\gamma_1^2-\frac{1}{48}\pi^4\gamma_1^2+\gamma_1^4+2\gamma\gamma_2-\frac{1}{48}\gamma\gamma_2\pi^4+\\
&-2\gamma\gamma_1^2\gamma_2+\frac{1}{2}\gamma^2\gamma_2^2-\gamma_1\gamma_2^2+\frac{2}{3}\gamma_3-\frac{1}{144}\pi^4\gamma_3+\frac{2}{3}\gamma^2\gamma_1\gamma_3+\frac{2}{3}\gamma_1^2\gamma_3+\\
& +\frac{2}{3}\gamma\gamma_2\gamma_3+\frac{1}{6}\gamma_3^2-\frac{1}{6}\gamma^3\gamma_4-\frac{1}{3}\gamma\gamma_1\gamma_4-\frac{1}{12}\gamma_2\gamma_4-\frac{1}{30}\gamma^2\gamma_5-\frac{1}{180}\gamma\gamma_6+\\
&-\frac{1}{1260}\gamma_7+\gamma^4-\frac{\pi^4}{96}-\frac{1}{96}\pi^4\gamma^4+\frac{23}{215040}\pi^8\Bigg)^{-\frac{1}{4}}-\frac{1}{4}\Bigg]^{\frac{1}{2}}
\end{aligned}
\end{equation}

\begin{equation}\label{eq:20}
t_1\approx 12.161258748655529488677538477512\ldots. \nonumber
\end{equation}
\end{document}